\documentclass[12pt]{article}
\usepackage{amsmath}
\usepackage{amssymb}
\usepackage{amsthm}
\usepackage{amscd}
\usepackage{amsfonts}
\usepackage{dsfont}
\usepackage{cases}
\usepackage[applemac]{inputenc}
\usepackage{enumerate}
\usepackage[all]{xy}
\usepackage{esint}
\usepackage{graphicx}
\usepackage{url}
\theoremstyle{plain}
\newtheorem{theor}{Theorem}{\normalfont}

\theoremstyle{definition}

\newcommand{\N}{\mathbb N}

\newcommand{\Z}{\mathbb Z}
\newcommand{\E}{\mathrm E}

\newcommand{\M}{\cal M}
\newcommand{\A}{{\cal A}}

\newcommand{\Om}{\Omega}

\begin{document}
\title{Galton-Watson processes and their role as building blocks  for branching processes }
\date{}
\author{F. Thomas Bruss\\Universit\'e Libre de Bruxelles}
\maketitle

\begin{abstract} This article is an essay, both expository and argumentative, on the Galton-Watson process as a tool in the domain of Branching Processes.  It is at the same time the author's ways to honour two distinguished scientists in this domain, both from the Russian Academy of Science, and to congratulate them for their special birthdays coming up very soon.

The thread of the article is the role, which the Galton-Watson process had played in the author's own research.  
We start with article on a controlled  Galton-Watson process. Then we pass to  {\it random absorbing processes}, and also recall and discuss a problem in medicine. Further questions will bring us via the Borel-Cantelli Lemma
to $\varphi$-{branching processes} and extensions. To gain more generality we then look at bisexual Galton-Watson processes. Finally we briefly discuss relatively complicated {\it resource dependent branching processes} to show that, here again, using Galton-Watson reproduction schemes (whenever reasonable) can be a convincing approach to new processes which are then sufficiently tractable to obtain results of interest.


\smallskip\noindent
{\bf Math. Subj. Index:} Primary 69J85, 60J05; secondary 60G40.
\end{abstract}

\theoremstyle{plain}

\section{Introduction} Would you remember the name of the first author you cited in your very first publication? The author did not have to look it up. It was Zubkov (1970) in his first publication (B.\footnote{In this article the author's name will be abbreviated throughout by B.} (1978)), and the same was  true for his thesis
a year earlier.

\smallskip
Among many distinguished Russian  authors in the field of branching processes, two will be honoured in this essay, and this  for two reasons. First, both have special birthdays coming up.  One of them is, as the reader guessed already,  Andrei  M. Zubkov (born December 30, 1946), and the second one is Vladimir A. Vatutin (January 1, 1952). Both are from the Russian Academy of Sciences in Moscow. Second, the author sees them both as masters in contributing to, and dealing with, the Galton-Watson process
and its extensions.
\subsection{The Galton-Watson process}

 Let $(Z_n)$ be a Galton-Watson process (GWP) with $Z_0=1$ and reproduction mean $m= \sum_{k=1}^\infty k\, p_k,$  where $p_k$ denotes the probability, that a particle will have exactly $k$ offspring in the next generation. We suppose $p_0>0$ and $p_0+p_1<1$ so that the limiting behavious of $(Z_n)$ is non-trivial.
 
 Since $\{0\}$ is an absorbing state we know that $P(\lim_{n \to \infty}Z_n=0)=\lim_nP(Z_n=0).$ Hence $q:= \lim_nP(Z_n=0)$ gives the probability of final extinction of the GWP.  Moreover we know that, due to  the independence of reproduction of particles in all generations, 
\begin{align} \forall j \in \N:  P(\lim_{n \to \infty}Z_n=0| Z_j=k)=q^k ,~k=1, 2, \cdots \end{align}
which, if $q<1,$  tends exponentially quickly to zero as k increases.  This implies that if $(Z_n)$ dies out, then it is very probable that this happens early. Time itself is irrelevant in the sense that, as we have just recalled,  $P( Z_n=0| Z_j) \to q^k$ for all intermediate times $j.$ Hence, the smaller the number $k$ of living  particles in any generation, the larger is the probability that the process will finally die out.

Consequently, a natural question is whether one can find a time-dependent ``minimal" cutting rule  (i.e. a rule to take off in each generation a certain number of particles) such that, with supercritical reproduction,  the final extinction would still be almost sure?

If $m\le 1$ then we know that $Z_n\to 0~a.s.$ so that no cutting is needed at all. This confines our interest to the case $m>1.$
Also, according to (1) we feel immediately that additive constants should be irrelevant for this question, and that therefore  ``minimal" should mean in fact ``of minimal order".
\smallskip

What then is this minimal order?
If the author guesses right, this must have been, in its essence,  the problem suggested to Zubkov by his prolific academic father Sevast'yanov.  

\section{Zubkov's model and modifications}

Zubkov (1970) made this question precise by proposing the following model:

\smallskip
\noindent Let $N(z)$ denote the random number of offspring generated in one generation by $z$ particles which all reproduce independently according to the same law $\{p_k\}_{k=0, 1, 2, \cdots}$ 
Let $(Z_n)$ be the corresponding  GWP with $Z_0=1.$ Further, let  \[Z_n:= N(Z_{n-1}),~n=1, 2, \cdots\] be the short notation for saying that, given $Z_{n-1}=z,$ the random variable $Z_n$ is distributed like $N(z).$  Then

\begin{theor} \normalfont
(Zubkov 1970) 
\theoremstyle{plain} Let $g$ be an integer-valued deterministic function from $\{0, 1, 2, \cdots\}$ into $\{0,1, 2, \cdots\}$ with $g(0)\ge 1$, and let \begin{align}\tilde Z_n:=\begin{cases}\tilde Z_0:=Z_0=1&, \text{if~}n=0\\\min\{g(n), N(\tilde Z_{n-1})\}&, \text {if~} n=1, 2, \cdots, \end{cases}\end{align} where $N(x)$
denotes the random number of offspring generated by $x$ particles under independent reproduction according to the  original GWP-law $\{p_k\}.$ Then 
\begin{align} \tilde q:= P(\lim_n \tilde Z_n=0) =1 \iff \sum_{n=1}^\infty q^{g(n)} = \infty.\end{align}

\end{theor}
 Zubkov's answer was thus as crisp and clear as one could possibly hope for! His result was an if-and-only-if divergence condition, linking in (3)  the upper deterministic truncation function $g(\cdot)$ of the modified (truncated) process with the extinction probability $q$ of the original GWP. 
  \subsection{The expected value of the growth curve}
  
  For the author, Zubkov's result was  ``before his time", since, before he had seen it,
he had been looking independently at the following model: 

\smallskip

For any non-negative integer $\ell,$ let $A_n(\ell), n=1, 2, \cdots$ be {\it arbitrary} random variables taking value in $\{0, 1, \cdots,\ell\}.$  The author thought of $A_n(\ell)$ as a random number of {\it absorbed} particles among $\ell$ and called the process $(Z_n^A)_{n=1, 2, \cdots}$ defined by the recursion
\begin{align} 
Z_n^A =\begin{cases}
Z_0=1&,  n=0\\ N(Z_{n-1}^A)-A_n\big(N(Z_{n-1}^A)\big) &, n = 1, 2, \cdots
\end{cases}
\end{align}
 {\it branching process with a random absorbing process}. 
 
For the absorbed number $A_n,$ the intonation is here on {\it arbitrary} within the set of possible values. In particular we allow the $A_n(\ell)$ not only to depend on $n$ but also on the whole history of the process $Z_1^A, Z_2^A, \cdots, Z_{n-1}^A.$
The definition of the absorbing process $A_n(\cdot)$ is thus not specified. Hence Zubkov's controlled branching process is the special case $$A_n\big(N(Z_{n-1}^A)\big)=\big[ g(n)-
N(Z_{n-1}^A)\big]^+,$$namely the number of particles in generation $n$, if it exceeds $g(n),$  will be truncated to the number $g(n)$ of particles which may reproduce. 
As the author thinks, because of (3), this is an important special case indeed. 

\smallskip
Here now a related result:
\begin{theor} \normalfont
(B. 1978) 
\theoremstyle{plain}
Let $g: \N \to [1,\infty]$ and, as before, $q$ denote the extinction probability of the original (unmodified)
GWP $(Z_n).$ Further suppose that $\E(Z_n^A|Z_n^A>0)\le g(n)$ for almost all $n.$ Then we have 
\begin{align}\sum_{n=1}^\infty p^{g(n)} = \infty {\rm~for ~some~} p\in ]0,q]~\implies~q_A:= P(\lim_{n\to \infty} Z_n^A=0)=1. \end{align}
\end{theor}

 We note that the lhs-sum in (5) is increasing in $p$ so that (in order to display the analogous part in Zubkov's theorem) we can replace $p$ in the lhs sum directly by $q.$ Alternatively, if we replace  in (5) the word {\it some}  on the lhs by {\it any},  then it is easy to show that the implication holds also under the weaker boundary condition $\E(Z_n^A)\le g(n)$ for almost all $n.$ Depending on the situation, both forms can be of interest for possible applications, as explained below.
 
Moreover, a moment's reflection shows that we cannot hope for an equivalence relation as seen in Zubkov's theorem. (Theorem 1.) 

So, for instance, we may assume that the absorbing process $(A_n)$ is such that it either absorbes all particles at once (disaster), or alternatively, it absorbs no particle at all. Then $\E(Z_n^A| Z_n^A>0)$ grows exponentially quickly since $m>1.$ However, if a disaster will arrive at some time almost surely, then clearly $q_A=1.$ To obtain the latter it suffices (for instance)  to suppose that a first disaster arrives at time $k$, independently of the history up to time $k-1$, with probability $\delta_k,$ and that $\sum_{k=1}^\infty \delta_k= \infty.$

\smallskip
We also  mention here that Schuh (1976) studied  an interesting related problem of a Galton-Watson which is absorbed as soon as the number of effectives drops below some lower bound.  Thus the idea in Schuh's setting  is to model a population which cannot survive unless it grows sufficiently quickly.

\subsection {Motivation}
\smallskip
The author's motivation behind Theorem 2 came from thinking of medical doctors, or medical staff, fighting against  bacteria or malignant cells in patients (B. (1979)). Medication may help to reduce or stop their growth. 

Often enough however, medication has bad side effects on the patients and also reduce their quality of life. Hence one would like to keep the dosage of the medication always as weak as possible, without losing sight of the objective to stop or control growth of malignant cells or bacteria. Typically, doctors would not know a deterministic control function $g(n)$, but they may think of repetitive experiments under the same circumstances, giving information on $\E(Z_n^A)$ or $\E(Z_n^A|Z_n^A>0),$ ideally in laboratory conditions not endangering the patients. Here the first option allows a quicker collection of data since observations become superfluous for extinct growth trajectories, but the second one is, in several aspects,  more informative. 

Clearly, we all see that the objective of finding a minimal cutting rule is not well-defined in the sense
that there exists no most slowly diverging series. Hence there can be no optimal control as such. However, this truly lies in the nature of the problem, and the essence of the problem remains of interest despite this intrinsic formal weakness of the problem itself.  

\smallskip

To exemplify this line of thought in another context, we refer to Dietz (1973), of whom the author had heard not long before he submitted his thesis
in 1977. Indeed, Dietz, using a different model, gave a convincing description of a real-world problem of this kind. Dietz studied the effect of the so-called ``sterile male technique", in which the growth of supercritical populations of insects is controlled by releasing in certain time intervals infertile male insects. The objective was to control the growth of the population of insects, not to extinguish the whole population. (Dietz became the main external examiner on the board of examiners at the author's thesis defence.)

\smallskip
For real-world applications of such a type, it will likely be necessary to think about estimation procedures for unknown parameters. Here control is usually the result of observations of the population's growing behaviour (observed in the growth trajectories), and the target is to understand the essential parameters which seem to govern these trajectories. In practice, it is often difficult, sometimes even impossible, to observe trajectories over a long time. Several studies have been made in the
direction to overcome this problem. See e.g. the more recent research on  stepwise estimation methods based on progenitors of each generation (Gonz\'ales et al. (2016)).

\section {Criticality}
We also note that questions of control of this type can be translated into questions of criticality. For a branching process, in which $0$ remains an absorbing state, the question becomes:  When does a process linger long enough around some boundary function (and which one?) so that fluctuations will bring the process finally down into the absorbing state $0,$ or, alternatively allow for unlimited growth? 

Questions of this type must have been played in the 1970-ties and 1980-ties an important role for the motivation  of the work of Vatutin and Zubkov whom we honour in this article. See
Vatutin (1977) and all the work of Vatutin and Zubkov cited in the References.  Again, some of the work in their collaboration was, as the author understands, suggested and/or initiated by Sevast'yanov. 

The collaboration of Vatutin and Zubkov over many years was very successful and had seemingly a solid basis of  mutual respect. When reading the author's comment
that Zubkov's Theorem 1 was {\it before his time}, Vatutin (private communication) kindly added, that this is in his opinion also true for Zubkov's work on limit distributions of the distance to the nearest mutual ancestor  (Zubkov (1975)). The author was not aware of Zubkov's work on this subject. However, he also has heard repeatedly, that questions concerning the nearest mutual ancestor have attracted and do attract a great deal of interest.

\smallskip
Returning to the issue  of fluctuations, our question above creates  also an interesting  link with the later work of Afanasyev et al. (2005). 

These authors studied branching processes under assumptions flowing out of Spitzer's condition in fluctuation theory for random walks. Their processes are more specific than  branching processes submitted to arbitrary absorbing processes, as we have suggested in B. (1978), and which have some advantages through their generality.  However, the results of Afanasyev et al. are stronger and more useful in other cases. 
This strength lies in a conditional functional limit theorem which these authors obtained by studying the behaviour of the survival probability, resulting in solid conditional functional limit theorems for the generation size process.
\section{$\varphi$-branching processes}
Returning again to earlier work, mainly of Vatutin and Zubkov, it is clear that the wish to understand the evolution of populations in which the population size at each generation needs to be controlled is behind the objective of creating suitable models. This wish may have motivated Sevast'yanov and Zubkov (1974) to create their so-called $\varphi$-branching processes (equally known as $\phi$-branching processes with capital Phi) and to study them thoroughly. Vatutin and Zubkov, and seemingly many other authors (as e.g. Dietz (1973) and several other names cited in this essay, including the author), must have been guided by this wish to obtain a control which is hoped to be  beneficial in some sense.

\smallskip We briefly recall the definition.
 
 A  (discrete-time) $\varphi$-branching process is a non-negative  integer-valued stochastic process $(\mu(t))_{t=0, 1, 2, \cdots}$ which can be defined on the basis of a double-array of i.i.d. non-negative integer random variables $\{\xi_k^t\}_{k,t \in \N}$ and a deterministic function $\varphi: \N \to \N. $ The process is then determined recursively by $\mu(0)=1$ and
 \begin{align} 
 \mu(t)= \sum_{k=1}^{\varphi(\mu(t-1))}~\xi_k^{t-1},~{\rm if~} t=1,2, \cdots.
\end{align}
 Here each $\xi_k^{t-1}$ is interpreted as the number of offspring of the $k$th particle in the $(t-1)$th generation.
 
 Most of the $\xi_k^{t-1}$ of the double-array will thus turn out redundant in this definition of $\mu(t)$ (unless the image of $\varphi$ is infinite). Moreover, clearly, if $\varphi$ is the identity function, the definition recaptures the definition of the GWP. Similarly, we see that immigration and emigration, and mixtures of these, can now be modelled by the corresponding choice of the function $\varphi.$
 
\smallskip

This model and the results of Zubkov and Vatutin seemingly motivated Yanev (1976)  to study his interesting generalisation of $\varphi$-processes for {\it random} control functions $\varphi.$ This again stimulated the author (B. 1980) to try to find  an extension of Yanev's version where the law of the control function $\varphi$ need not necessarily be specified as such, but only in an ``expectational sense". Recall that B.'s model of 1978 cannot be retrieved here
since there the $A_n(\ell)$ may depend on the whole history of the growth curve.

The extension of Yanev's version in the direction of ``expectational information'' required more, however, and this rejoins our question stated in the first paragraph. Starting from {\it expectations}, what can we do in order  to deal with the accumulated probabilities for a process to get absorbed in $0$, 
and then to estimate them?  The author could finally provide an answer for the first question relating to Yanev's model, which he thought to be useful. It was based on his modification of the Borel-Cantelli lemma, which we recall below.  

\smallskip The lemma of B. (1980)) is as follows:  
\section{A counterpart of the Borel-Cantelli Lemma}
Let $(E_n)_{n=1, 2, \cdots}$ be a sequence of events defined on some probability space $(\Om, \A, P)$ and, as usual, let $\{E_n {\rm~ i.o. }\}$ denote the event that infinitely many $E_n$ occur. The first part of the Borel-Cantelli Lemma (BCL) says that the divergence of the series $\sum_nP(E_n)$ is a necessary condition for $P(E_n~i.o)=1$  to hold. If one speaks of the second part of the BCL, one typically thinks of  the sufficient condition \begin{align} E_n \in \A {\rm~independent~for~all}~n,~{\rm and~}\sum_n^\infty P(E_n)=\infty \implies P(E_n {\rm ~i.o.})=1.~\end{align} It is well-known that the implication in (7) holds for more general conditions on the events $E_j$, such as for instance pairwise independence or some other weaker forms of independence (see e.g. Bir\'o and Curbelo (2020) for  a recent review).

However, for the mentioned objective of generalising the results of Yanev (1976) in an expectational form without specifying the law of a random control $\varphi$,  any independence assumption seems difficult to defend as realistic in the view of possible applications.
This lead to the idea to concentrate on the one property which is true for all branching processes without immigration,
namely that the state $0$ (extinction) is absorbing. Hence, if we let  $E_n:=\{Z_n^A=0\}$, we have for the process $(Z_n^A)$ together with the definition 
$$ \forall \,n, s \in \N: E_n ~{\rm implies~} E_{n+s}$$ This means  $E_n\subset E_{n+1} \subset \cdots \subset E_\infty$, and thus \begin{align} P(E_n)\uparrow P(E_n~i.o.) = P(\lim Z_n^A =0).\end{align} We then have the following
version of the Borel-Cantelli Lemma

\begin{theor} \normalfont (B. 1980) Let $(E_n)$ be a sequence of events defined on a probability space $(\Om, \A, P)$, and let $\overline E$ denote the complement of the event $E.$ Then
$P(E_n~ i.o.)=1$ if and only if there exists a strictly increasing sequence $(t_k)_{k=1,2,\cdots}$, with $t_k \in \N,$ such that \begin{align} \label{B1980}\sum_{k =1}^\infty P\left(E_{t_k}\big|\,\overline{E_{t_{k-1}}}\,\right) = \infty.\end{align}
\end{theor}

The essence of the idea behind this Theorem is as follows. 
If we want to show, that infinitely many events of a certain type will occur, then it suffices of course to show that there is no last one of them. This means that after any chosen time $k$ another event will follow almost surely, and that just {\it one} will do. 
But then, why not suppose that all $E_n$ entail each other thereafter? 
Moreover, there is no other constraint on the $t_k$  than being in $\N,$ and $(t_k)$ being strictly increasing. Therefore we may ask, why not directly putting $t_k=k?$ This is de facto equivalent. Nevertheless, there is some benefit in the more general-looking  formulation. Sometimes the choice of the differences $t_{k+1}-t_k$ allows us to see more easily whether accumulated absorption probabilities build up sufficiently quickly to obtain a divergent series. 

\smallskip
When preparing this essay, the author discovered that Barndorff-Nielsen gave already in 1961 versions of the BCL which are related with Theorem 3. It is a pleasure to be able to add his reference now. Unlike Theorem 3,
Barndorff-Nielsen's results are in several Lemmata. All  sum terms in there are probabilities
of intersections. To compute or estimate these, one would write them as products of a conditional probability and an absolute probability, and since sums are easier to grasp than sums of products, Theorem 3 seems a priori preferable. Moreover, since the product of probabilities is bounded above by each factor, the newer single Theorem 3 is expected to be stronger.\footnote{This is no claim. The author will have to check.}

\smallskip Theorem 3 (Counterpart of the BCL) has proved to be useful not only in problems  motivated by  Branching processes, as explained before and explicitly done in Example  2. of B.(1980),
 but also in the context of rather different problems. For a recent application, we refer to  Makur et al. (2020) on a problem of broadcasting on a two-dimensional regular grid, and to their proof of  Proposition 4.  See also Feldman and Feldman (2020), Wirtz (2019), Cohen and Fedyashov (2018), Feldman and Souganidis (2017), Bertacchi et al. (2014) and others for further applications.

\section {Bisexual Galton-Watson processes} \label{Bisex} It is well-known, that the study of bisexual GWPs is in general more complicated than that of the (classical) GWP. The reason is that both sexes intervene, and that the so-called mating functions $\M$ defined on  the sexual behaviour of the two sexes can differ in many ways. Thus the law $\{p_k\}_{k=0,1,\cdots}$
and the mean reproduction rate of males  (respectively  females) are not enough to  yield extinction criteria in a necessary-and-sufficient form. 

However, if we consider now each couple brought together by the mating function as one {\it mating unit} then we obtain a tractable extinction criterion. As we will see, this criterion is strongly related with the extinction criterion for the  classical GWP.
Indeed, let $(X_n,Y_n)_{n=0,1,2 \cdots}$ be a bi-variate process, where $X_n~(Y_n)$ is seen as  the number of female particles (male particles) in the $n$th generation. We define the process $(Z_n)$ by
\begin{align} \label{Zn}Z_n:= {\M(}X_n, Y_n), ~ n=0, 1,2, \cdots \end{align}
where the so-called mating function ${\M} : \N^2 \to \N,~ (x,y)\to {\M}(x,y)\in \N,$ describes the number of mating units (mating pairs) $x$ females and $y$ males will form. Here we make the reasonable assumption that for each fixed $x$, ${\M}(x,Y_n)$ is non-decreasing in $Y_n,$ and similarly for each fix $y$, ${\M}(X_n,y)$ is no-decreasing in $X_n.$ Each mating unit is now supposed to reproduce independently according to the same GWP-law $\{p_k\}_{k=0,1,2, \cdots}$, where each offspring will be, independently of all other offsprings,  male  with probability $0<\alpha<1$ 
respectively female ith probability $1-\alpha.$ Offspring in the $(n+1)th$ generation are then $X_{n+1},Y_{n+1}$, and $Z_{n+1}$ assumes the value 
${\M}(X_{n+1},Y_{n+1}).$ This recursion determines the bisexual GWP (BGWP).

\smallskip
We then have:

\begin{theor}[B. 1984]\label{B. (1984)} Let the \emph{average reproduction rate} of $k$ mating units (measured in terms of mating units) be defined by
\begin{align} m(k)=\frac{1}{k}\E(Z_1|Z_0=k),~ k=1, 2, \cdots.\end{align} If $m(k)$ is bounded for all $k$, and $m(k)\le 1 $ for all $k$ sufficiently large, then \[Q=\lim_{n\to\infty} P(Z_n=0)=1.\]\end{theor}

 Daley (1968) had already studied the special cases ${\M}(x,y)=x \min(1,y)$ and ${\M}(x,y)=x \min(x, d y)$, where $d\in \{1,2, \cdots\},$ and in this case the above condition can be turned into an if-and-only-if condition.  Also, our preceding Theorem is partially covered 
 in the work Sevast'yanov and Zubkov (1974) if one supposes in addition $p_0>0,$ and moreover $m(k)<1$ for all $k$ sufficiently large.
 
 \smallskip
 As the reader will guess, one can  show (for example by using the counterpart of the BCL) that there can exist no if-and-only-if condition for extinction without restrictions on the mating function ${\M}$. By definition, this  BGWP cannot survive unless all $Z_n$ are strictly positive. Hence  it is worth noting that, as the preceding Theorem shows, the GWP-analogy holds  in great generality for the {\it necessary} condition for survival of the bisexual GWP.
 This is, among other contexts, important for the following Section.
 
 \section{Resource dependent branching processes}
 
 Resource dependent branching processes (RDBPs) are  branching processes which have been created with a relatively
 ambitious motivation, namely to study the development of human populations, i.e. populations  with typical human preferences. These are that individuals want to have an environment in which they and their children can survive safely and peacefully, and this with a standard of living as high as possible.
 
 These objectives are not always compatible which each other, and then, by definition, survival has the first priority. 
 It is the longer-term development of large poulations which will interest us. This is why, whenever a comparison with a result from a Galton-Watson-type process is used, we can, as shown in Section  \ref{Bisex}  neglect the feature of bisexual reproduction.
 
 Speaking now about RDBPs (first presented in rudimentary form at the 11th SPA-Conference in Clermont-Ferrand , see B. (1984)), we get further away from the work of Vatutin and Zubkov, and, as far the author knows, further away from the work of most other colleagues in the domain of Branching processes. In  this last section I try, now selfishly, to stimulate some interest in RDBPs. The difficulties in studying RDBPs stem from
 their particularities explained below, and any shared interest for the questions which turn up naturally, would be much appreciated.
 \subsection{Particularities of RDBPs}
In RDBPs, particles have to work in order to be able to live and to reproduce. They live ({\it consumption}) from resources they produce or which were left by their ancestors, and they create ({\it production}) new resources by their work, which goes then into a common resource space for all. 
Individual demands to obtain resources are called {\it claims.} These individual claims are submitted to the (current) society, which has its rules how to satisfy claims. 

Individuals have also the possibility to protest or interact, and thus a vote
on the society form which will be in command in the next generation. B. and Duerinckx (2015) specified the protest by assuming, that all individuals whose claims are not satisfied, prefer to either emigrate before having descendants, or,
with the same effect, to refuse to reproduce. However, the models are sufficiently flexible for other types of protest and control.

The idea is that the whole development of the population in time should be described by a sequence of such RDBP's. Hence RDBPs are themselves only the parts in a sequence determining the development of the population.  The parameters of RDBPs  may change, sometimes at each new generation. but the priority of wanting to survive ``forever" (in the sense of branching processes) is maintained.
Since the latter is a terminal event, a temporary society rule is in general not  important. To allow for a minimum of structure, it is ruled that
a society form under which, with its current rules and currently observed parameters, could not survive forever with a positive probability, must adapt accordingly to make this goal possible. This over-all priority is insuffient to describe the development of the population generation by generation, but it does allow to answer some questions of macro-economic interest.
 
 More details would exceed the scope of the present article, and the interested reader is referred to B. and Duerinckx (2015). Their main results of macro-economic are a conditional and unconditional form of the {\it Theorem of Envelopment.}
\subsection{Expected stopping times}
It is a chacteristic feature of such RDBPs that stopping times come in very naturally. Societies must stop distributing resources as soon as all those created by their ancestors are exhausted, or societies cannot increase productivity
 of their (learning) individuals above respective limits. RDBPs can thus be seen, on a larger scale, also as stochastic processes, or sets of stochastic processes, interacting by stopping times.
The goal is to see how societies develop under different society structures (rules imposing stopping times), how the outcome may influence the successive choices of society forms, and how the corresponding society survival probabilities behave.  
For these three types of questions the paper by B. and Duerinckx (2015) gives for a single population several answers which we found interesting, in particular the {\it Theorem of envelopment}\label{TE} for societies.

A major fourth objective looming behind all these questions is to better understand how a peaceful co-habitation of {\it sub-populations} with different rules
can work out. This question translates into how certain equilibria between sub-populations (in the number of effectives) can be reached.   Sub-populations may for instance include immigrants, which are typically quite different in their behaviour of production and consumption. 

Progress has been made into this direction by using a versatile upper bound for the expected value of a stopping time defined by random variables underlying different distributions. The inequality is as follows:

\smallskip
 \noindent Let $X_1, X_2, \cdots, X_n$ be positive random variables that are jointly continuously distributed and such that each $X_k$ has an absolutely continuous distribution function $F_k,$
 and let their increasing order statistics be denoted by $X_{1,n} \le X_{2,n} \le \cdots \le  X_{n,n}$.
 Let $s>0$ be a fixed real number, and let  the random variable $N(n,s)$ be defined \begin{align}
N(n,s)=\begin{cases} 0, {~\rm if}~ X_{1,n}>s,\\
\max\{k \in \N: X_{1,n}+X_{2,n}+ \cdots +X_{k,n} \le s\}, \,{\rm otherwise.}\end{cases}\end{align}
We note that $N(n,s)$ is a stopping time, and we have the following upper bound for its expectation (see Steele (2016), and for more details see B. (2021).)
\begin{align}\E\left(N(n,s)\right) \le \sum_{k=1}^n F_k(t),
\end{align}where
$t:=t(n,s)$ is the unique solution of the equation
$\sum_{k=1}^n \int_0^t x dF_k(x)= s.$

\smallskip
To see the link with what we explained before, the number $s$ is thought of as being a total available resource space which is used to satisfy the mentioned claims from the different sub-populations
by serving with priority the most modest claims first.
Hence if there are, for instance, three sub-populations
of respective sizes $n_1, n_2,$ and $n_3$ with $n=n_1+n_2+n_3$, and if claims have, within each sub-population, the same marginal distribution, then the rhs-sum in (13) can be replaced by $n_1F_1(t)+n_2F_2(t)+n_3F_3(t)$.  

The great advantage of this bound in (13) is that it holds {\it without any independence assumption}, so that it fits nicely with our general definition of a RDBP.
If one assumes moreover independence, a better bound would be available (see B.(2021) for details),
but this assumption is no  real-world assumption.

\smallskip
Concerning the fourth objective announced above, we only mention, that more general answers are still not yet at what one would really like to have. It would be nice to obtain rather general sufficient criteria for the existence and feasibility of computation of equilibria for the ratio of the number of effectives of the given sub-populations.

\subsection{Back to the Galton-Watson process}
Now, where exactly do Galton-Watson processes come in?
RDBPs, with their unavoidable stopping times, become  different from so many interesting branching processes we have cited, and also those of Haccou et al., Klebaner, and many other colleagues. Nevertheless, the Galton-Watson process is a distinguished building stone for more complicated branching processes.
We try to use its simple definition and its convenient extinction whenever we can. 

Clearly, no real-world RDBP can be Galton-Watson process. However, as long as interacting processes with their stopping time are not ``stopped", it makes sense to assume that reproduction is internally (i.e. within a given population of a certain type)
determined by i.i.d. random numbers of offspring. This opens, among other advantages, also the way to (stopped) martingales. Consequently, we can achieve to {\it quantify} intermediate steps and results, rather than having to be satisfied with
qualitative statements, since it is the quantified statements which typically give the stronger results. 

This is often much easier for each GWP in the setting than for other building blocks. 
See for example the proofs of the two versions of the mentioned Theorem of envelopment (for one population) in B. and Duerinckx (2015).

\smallskip
As an example to see how this can also work with interaction between building blocks, and specifically for an ``outer'' block and an ``inner'' block, recall again Section 6 on bisexual reproduction. There we had subordinated (or inserted) the $X_n$ males and $Y_n$ females into a process of mating units $Z_n,$ and obtained a GWP-like extinction criterion for the process of mating units $(Z_n)$.

\smallskip Please note that neither $(X_n)$, nor $(Y_n)$ nor $(Z_n)$ are Galton-Watson processes on their own. 
This was not really required. If $(Z_n)$ survives, then this means by definition for the real-world that both $(X_n)$ and $(Y_n)$ survive. If they do, they will behave like GWPs with limiting ratios $\alpha/(1-\alpha)$ because of the independent sex-selection assumption for the newborns. If not, all get extinct anyway. We did not need to specify the mating function $\M$ to get to the essentials.
\section{ Conclusion}

In conclusion, with all what we know about Galton-Watson processes, extinction criteria, asymptotic growth, their relationship with martingales, it seems good advice to always think of them, a little bit in the way computer scientist think of black boxes. Perhaps we should use them as building blocks whenever this is, in their right place, compatible with a given problem.

\smallskip
Recalling the special anniversaries of Zubkov and Vatutin, we end with a reference to music. Mozart has been cited to have said {\it Viel herrliche Musik kann noch in C-dur geschrieben werden.} (A  lot of beautiful music can still be written in C-major.) Whether  or not it was Mozart who said it first is not relevant. No doubt, it is true. 
GWPs cannot be for branching processes exactly what C-major is for music. Still, feelings do not come only out of the blue.

\bigskip{\bf Acknowledgement}

\medskip
\noindent The author has now heard  that the translation of the article Zubkov (1970), which H.-J. Schuh showed him in 1976/7, was owed to F. Klebaner. Thank you both.

\bigskip
\smallskip{\bf {References}}
~
\smallskip
 
Afanasyev V., Geiger J., Kersting G., Vatutin V. A. (2005)
`` Criticality for branching processes in random environment'',
Ann. Probab. Vol. 33 (2): 645-673

\smallskip Barndorff-Nielsen O. (1961) ~``On the rate of growth of partial extrema of independent, identically distributed random variables'', Amanuensis of Mathematics, Techn. Scient. Note 5, 1-17.

\smallskip
Bertacchi D., Machado F. P. and Zucca F. (2014)~``Local and global survival for non-homogeneous random walk systems on $\Z$'', Adv. in Appl. Prob., Vol. 46, 256-278.

\smallskip Bir\'o C. and Curbelo I. R. (2020) ~``Weak independence of events and the converse of the Borel-Cantelli Lemma", Expositiones Mathematicae, published online; available at arxiv: 2004.11324v2

\smallskip Bruss F. T. (1978)~``Branching Processes with Random Absorbing Processes", J. Appl. Prob.,  Vol. 15,
54-64.

\smallskip Bruss F. T. (1979)~``How to Apply a Medicament, if its Direct Efficiency is Unknown'',
Annales de la Soci\'et\'e Scientifique de Bruxelles, Vol.  93, No 1,  39-54.

\smallskip Bruss F. T. (1980)~``A Counterpart of the Borel-Cantelli Lemma'', J. Appl. Prob., Vol. 17,  1094-1101.

\smallskip Bruss F. T. (1984)~``A Note on Extinction Criteria for Bisexual Galton Watson Processes'',
J. Appl. Prob., Vol. 21, 915-919.

\smallskip

 Bruss F. T. (1984b) Resource Dependent Branching Processes, 11th Conf. on Stoch. Proc.
Applic. 1982, Abstract, Stoch. Proc. Applic., Vol. 16 (1), p. 36.

\smallskip  Bruss F. T. and Duerinckx M. (2015) ~``Resource Dependent Branching Processes and the Envelope of Societies'', Annals of Appl. Probab.,  Vol. 25 , No 1, 324 - 372.

\smallskip Bruss F. T. (2021)~``The BRS-inequality and its applications", Probability Surveys, Vol. 18, 44-76.

\smallskip Cohen S. N. and Fedyashov  V. (2018), ``Ergodic BSDEs with jumps and time dependence''
available at arXiv:1406.4329v2 

\smallskip  Daley D. J. (1968), `` Extinction conditions for certain bisexual Galton-Watson branching processes'', Z. Wahrscheinl. und Verw. Geb., Vol. 9  (4),  315–322. 
 
 \smallskip Dietz K. (1973),~`` Simulation Models for Genetic Control Alternatives'', Lecture Notes in Biomathematics,
Vol. 5, 299-317.

\smallskip Feldman N. D. and  Feldman O. N. (2020) ~``Convergence of the quantile admission process with veto power'',
Stoch. Proc. and Applic., Vol. 130 (7), 4294-4325.

\smallskip
Feldman W. M, and Souganidis P.E. (2017)
``Homogenization and non-homoge-nization of certain non-convex Hamilton-Jacobi equations",
Journal de Math\'e-matiques Pures et Appliqu\'es,
Vol. 108 (5), 751-782.
 
  \smallskip 
 Gonz\'alez M., del Puerto I., and Yanev G., (2017) ~
`` Controlled Branching Processes'', J. Wiley \& Sons, Vol. 2.

\smallskip Gonz\'ales M., Minuesa C., and del Puerto I.  (2016)~``Maximum likelihood estimation and expectation: maximization algorithm for controlled branching processes", Computational Statistics \& Data Analysis.,
Vol. 93, 209-227.

\smallskip 
Haccou P.,  Jagers P. and  Vatutin V.A. (2005) ~{\it Branching Processes - Variation, Growth and Extinction of Populations}, Cambridge University Press.

\smallskip
Makur A.,  Mossel E., and Polyanskiy Y. (2020)~``Broadcasting on two-dimensional regular grids'', published online, available at
arxiv: 2010.013 90v1.

\smallskip Schuh H.-J. (1976)~ ``A condition for the extinction of a branching process with an absorbing lower boundary'', J. Math. Biol. Vol. 3, 271-287.

\smallskip Sevast`yanov B. A. and Zubkov A. M. (1974)~``Controlled Branching Processes'', Theory Prob. Applic.,  Vol. 19, 14-24.

\smallskip
Steele J. M. (2016),~``The Bruss-Robertson Inequality: Elaborations, Extensions, and Applications",  Math. Applicanda, Vol. 44,  No 1, 3-16.

\smallskip 
Terelius H.  and  Johansson K. J. (2018)
~``Peer-to-peer gradient topologies in networks with churn'',
Journal
IEEE Transactions on Control of Network Systems
Vol. 5 (4), 
2085-2095.

\smallskip
Vatutin V. A. (1977),
``Asymptotic behaviour of the non-extinction probability for a critical branching process''
Theory of Probab. Appl., Vol. 22 (1), 140-146. 

\smallskip 
Vatutin V. A. and Zubkov A. M. (1985),
``Branching Processes I'', Itogi Nauki i Tekhniki Ser. Teor. Veronatn. Mat. Sat. Teor. Kibern. Vol. 23, 3-67.

\smallskip 
Vatutin V. A. and Zubkov A .M. (1993)
``Branching processes I'',  Journ.  Math. Sci., Vol. 67, 3407-3485.

\smallskip
Wirtz J. M. (2019)~``Coalescent Theory and Yule Tress in Time and Space'', Dissertation, Universit\"at zu K\"oln.

\smallskip 
Yanev N. M. (1976) ~``Conditions for degeneracy of $\phi$-branching processes with random $\phi$", Theory Prob.  Appl.,  Vol. 20, 421-427.

\smallskip 
Zubkov A.M. (1970) ~``A degeneracy condition for a bounded branching process'' (in Russian) Mat. Zametki Vol. 8, 9-14. English translation: Math. Notes {\bf 8}, 472-477.

\smallskip 
Zubkov A.M. (1975),
``Analogues between Galton-Watson processes and ${\phi}$-branching processes'',
Theory  Probab.   Appl., 19 (2), 309-331. 

\smallskip 
Zubkov A.M. (1976), ``Limit distributions of distance to the nearest mutual ancestor'', 
Theory  Probab.  Appl., 20 (2) 602-612.

\bigskip\medskip

\noindent Author's address

\smallskip

F. Thomas Bruss

Universit\'e Libre de Bruxelles

Facult\'e des sciences

Campus Plaine CP 210

B-1050 Bruxelles, Belgique

\smallskip
Thomas.Bruss@ulb.be

++32-2-650 5893
(++32-472-448801)

\end{document}